\documentclass{amsproc}
\usepackage{amssymb}

\usepackage{graphicx}
\usepackage{amscd}

\theoremstyle{plain}

\newtheorem{corollary}{Corollary}
\newtheorem{definition}{Definition}
\newtheorem{example}{Example}
\newtheorem{lemma}{Lemma}
\newtheorem{problem}{Problem}

\newtheorem{remark}{Remark}
\newtheorem{theorem}{Theorem}

\begin{document}
\title[On Banach spaces with the Tsirelson property]{On Banach spaces with the Tsirelson property}
\author{E.V. Tokarev}
\address{B.E. Ukrecolan, 33-81, Iskrinskaya str., 61005, Kharkiv-5, Ukraine}
\email{tokarev@univer.kharkov.ua}
\dedicatory{Dedicated to the memory of S. Banach.}

\begin{abstract}
A Banach space $X$ is said to have the Tsirelson property if it does not
contain subspaces that are isomorphic to $l_{p}$ $(1\leq p<\infty )$ or $%
c_{0}$. The article contains a quite simple method to producing Banach
spaces with the Tsirelson property.
\end{abstract}

\maketitle

A Banach space $X$ is said to have the Tsirelson property if it does not
contain subspaces that are isomorphic to $l_{p}$ $(1\leq p<\infty )$ or $%
c_{0}$.

The first example of a Banach space with such property was constructed by
B.S. Tsirelson [1].

The article contains a quite simple method to producing Banach spaces with
the Tsirelson property. Results were communicated at the Ukrainian
Mathematical Congress (August 2001, Kiev) and was announced in the Book of
Abstracts of the International Conference on Functional Analysis that take
place as a part of the Congress [2].

\section{Definitions and notations}

Let $\mathcal{B}$ be a (proper) class of all Banach spaces.

\begin{definition}
Let $X$, $Y\in \mathcal{B}$. $X$ is \textit{finitely representable} in $Y$
(in symbols: $X<_{f}Y$) if for each $\varepsilon >0$ and for every finite
dimensional subspace $A$ of $X$ there exists a subspace $B$ of $Y$ and an
isomorphism $u:A\rightarrow B$ such that 
\begin{equation*}
\left\| u\right\| \left\| u^{-1}\right\| \leq 1+\varepsilon .
\end{equation*}

Spaces $X$ and $Y$ are said to be finitely equivalent if $X<_{f}Y$ and $%
Y<_{f}X$.

Any Banach space $X$ generates classes 
\begin{equation*}
X^{f}=\{Y\in \mathcal{B}:X\sim _{f}Y\}\text{ \ and \ }X^{<f}=\{Y\in \mathcal{%
B}:Y<_{f}X\}
\end{equation*}
\end{definition}

For any two Banach spaces $X$, $Y$ their \textit{Banach-Mazur distance }is
given by 
\begin{equation*}
d(X,Y)=\inf \{\left\| u\right\| \left\| u^{-1}\right\| :u:X\rightarrow Y\},
\end{equation*}
where $u$ runs all isomorphisms between $X$ and $Y$ and is assumed, as
usual, that $\inf \varnothing =\infty $.

It is well known that $\log d(X,Y)$ defines a metric on each class of
isomorphic Banach spaces. A set $\frak{M}_{n}$ of all $n$-dimensional Banach
spaces, equipped with this metric, is a compact metric space that is called 
\textit{the Minkowski compact} $\frak{M}_{n}$.

A disjoint union $\cup \{\frak{M}_{n}:n<\infty \}=\frak{M}$ is a separable
metric space, which is called the \textit{Minkowski space}.

Consider a Banach space $X$. Let $H\left( X\right) $ be a set of all its
different finite dimensional subspaces (isometric finite dimensional
subspaces of $X$ in $H\left( X\right) $ are identified). Thus, $H\left(
X\right) $ may be regarded as a subset of $\frak{M}$, equipped with the
restriction of the metric topology of $\frak{M}$.

Of course, $H\left( X\right) $ need not to be a closed subset of $\frak{M}$.
Its closure in $\frak{M}$ will be denoted $\overline{H\left( X\right) }$.
From definitions follows that $X<_{f}Y$ if and only if $\overline{H\left(
X\right) }\subseteq \overline{H\left( Y\right) }$. Spaces $X$ and $Y$ are 
\textit{finitely equivalent }(in symbols: $X\sim _{f}Y$) if simultaneously $%
X<_{f}Y$ and $Y<_{f}X$. Therefore, $X\sim _{f}Y$ if and only if $\overline{%
H\left( X\right) }=\overline{H\left( Y\right) }$.

There is a one to one correspondence between classes of finite equivalence $%
X^{f}=\{Y\in \mathcal{B}:X\sim _{f}Y\}$ and closed subsets of $\frak{M}$ of
kind $\overline{H\left( X\right) }$.

Indeed, all spaces $Y$ from $X^{f}$ have the same set $\overline{H\left(
X\right) }$. This set, uniquely determined by $X$ (or, equivalently, by $%
X^{f}$), will be denoted by $\frak{M}(X^{f})$ and will be referred as to 
\textit{the Minkowski's base of the class} $X^{f}$.

\begin{definition}
For a Banach space $X$ its $l_{p}$-\textit{spectrum }$S(X)$ is given by 
\begin{equation*}
S(X)=\{p\in \lbrack 0,\infty ]:l_{p}<_{f}X\}.
\end{equation*}
\end{definition}

Certainly, if $X\sim_{f}Y$ then $S(X)=S(Y)$. Thus, the $l_{p}$-spectrum $%
S(X) $ may be regarded as a property of the whole class $X^{f}$. So,
notations like $S(X^{f})$ are of obvious meaning.

\begin{definition}
Let $X$ be a Banach space. It is called:
\end{definition}

\begin{itemize}
\item  $c$-\textit{convex,} if $\infty \notin S(X)$;

\item  $B$-\textit{convex,} if $1\notin S\left( X\right) $;

\item  \textit{Finite universal,} if $\infty \in S(X)$.
\end{itemize}

Let $\{X_{i}:i\in I\}$ be a collection of Banach spaces. A space 
\begin{equation*}
l_{p}\left( X_{i},I\right) =\left( \sum \oplus \{X_{i}:i\in I\}\right) _{p}
\end{equation*}
is a Banach space of all families $\{x_{i}\in X_{i}:i\in I\}=\frak{x}$, with
a finite norm 
\begin{equation*}
\left\| \frak{x}\right\| _{p}=(\sum \{\left\| x_{i}\right\|
_{X_{i}}^{p}:i\in I\})^{1/p}.
\end{equation*}

If $I=\mathbb{N}$, then instead $l_{p}\left( X_{i},\mathbb{N}\right) $ it
will be written $l_{p}(X_{i})$. If all $X_{i}$'s are equal to a given Banach
space $X$, then the notation $l_{p}\left( X\right) $ is used.

If $I$ consists of two elements only (say, $I=\{1,2\}$), then $l_{p}\left(
X_{i},I\right) $ is denoted by $X_{1}\oplus _{p}X_{2}$.

\begin{definition}
Let $I$ be a set; $D$ be an ultrafilter over $I$; $\{X_{i}:i\in I\}$ be a
family of Banach spaces. An \textit{ultraproduct }$(X_{i})_{D}$ is a
quotient space 
\begin{equation*}
(X)_{D}=l_{\infty }\left( X_{i},I\right) /N\left( X_{i},D\right) ,
\end{equation*}
where $l_{\infty }\left( X_{i},I\right) $ is a Banach space of all families $%
\frak{x}=\{x_{i}\in X_{i}:i\in I\}$, for which 
\begin{equation*}
\left\| \frak{x}\right\| =\sup \{\left\| x_{i}\right\| _{X_{i}}:i\in
I\}<\infty ;
\end{equation*}
$N\left( X_{i},D\right) $ is a subspace of $l_{\infty }\left( X_{i},I\right) 
$, which consists of such $\frak{x}$'s$\ $that 
\begin{equation*}
\lim_{D}\left\| x_{i}\right\| _{X_{i}}=0.
\end{equation*}
\end{definition}

If all $X_{i}$'s are all equal to a space $X\in \mathcal{B}$ then the
ultraproduct is said to be the \textit{ultrapower} and is denoted by $\left(
X\right) _{D}$.

An operator $d_{X}:X\rightarrow \left( X\right) _{D}$ that asserts to any $%
x\in X$ an element $\left( x\right) _{D}\in \left( X\right) _{D}$, which is
generated by a stationary family $\{x_{i}=x:i\in I\}$, is called the \textit{%
canonical embedding }of $X$ into its ultrapower $\left( X\right) _{D}$.

It is well-known that a Banach space $X$ is finitely representable in a
Banach space $Y$ if and only if there exists such ultrafilter $D$ (over $%
I=\cup D$) that $X$ is isometric to a subspace of the ultrapower $(Y)_{D}$.

In what follows a notion of inductive (or direct) limit will be used. Recall
a definition.

Let $\left\langle I,\ll\right\rangle $ be a \textit{partially ordered set}.
It said to be \textit{directed} (to the right hand) if for any $i,j\in I$
there exists $k\in I$ such that $i\ll k$ and $j\ll k$ .

Let $\left\{ X_{i}:i\in I\right\} $ be a set of Banach spaces that are
indexed by elements of an directed set $\left\langle I,\ll \right\rangle $.
Let $m_{i,j}:X_{i}\rightarrow X_{j}$ $(i\ll j)$ be isomorphic embeddings.

\begin{definition}
A system $\left\{ X_{i},m_{i,j}:i,j\in I;i\ll j\right\} $ is said to be an 
\textit{inductive (}or\textit{\ direct) system} if 
\begin{equation*}
m_{i,i}=Id_{X_{i}};\text{ }m_{i,k}=m_{j,k}\cdot m_{i,j}
\end{equation*}
for all $i\ll j\ll k$ ($Id_{Y}$ denotes the identical operator on $Y$).
\end{definition}

Let 
\begin{equation*}
X=\cup \left\{ X_{i}\times \left\{ i\right\} :i\in I\right\}
\end{equation*}

Elements of $X$ are pairs $(x,i),$ where $x\in X_{i}$. Let $=_{eq}$ be a
relation of equivalence of elements of $X$, which is given by the following
rule:

\begin{center}
$(x,i)=_{eq}(y,j)$ if $m_{i,k}x=m_{j,k}y$ for some $k\in I$.
\end{center}

A class of all elements of $X$ that are equivalent to a given $(x,i)$ is
denoted as 
\begin{equation*}
\lbrack x,i]=\{(y,j):(y,j)=_{eq}(x,i)\}.
\end{equation*}

A set of all equivalence classes $[x,i]$ is denoted $X_{\infty}$. Clearly, $%
X_{\infty}$ is a linear space. Let $\left\| [x,i]\right\| =\lim_{I}\left\|
m_{i,j}x\right\| _{X_{j}}$ be a semi-norm on $X_{\infty}$. Let 
\begin{equation*}
Null(X)=\{[x,i]:\left\| [x,i]\right\| =0\}.
\end{equation*}

\begin{definition}
A direct limit of the inductive system $\left\{ X_{i},m_{i,j}:i,j\in I;i\ll
j\right\} $ is a quotient space 
\begin{equation*}
\underset{\rightarrow }{\lim }X_{i}=X_{\infty }/Null(X).
\end{equation*}
\end{definition}

Let $X\in \mathcal{B}$; $\varkappa $ be an infinite cardinal; $%
dim(X)=\varkappa $; $\alpha $ be an infinite limit ordinal\ (so, $\alpha $
may be considered as a cardinal); $\alpha \leq \varkappa $.

\begin{definition}
An $\alpha $-sequence $\{x_{\beta }:\beta <\alpha \}$ of elements of $X$ is
said to be

\begin{itemize}
\item  Spreading, if for every $n<\omega $ ($\omega $ denotes the first
infinite cardinal or, equivalently, ordinal), every $\varepsilon >0$, every
set of scalars $\{a_{k}:k<n\}$ and any choosing of $i_{0}<i_{1}<...<i_{n-1}<%
\alpha $; $j_{0}<j_{1}<...<j_{n-1}<\alpha $ 
\begin{equation*}
\left\| \sum\nolimits_{k<n}a_{k}x_{i_{k}}\right\| \leq \left\|
\sum\nolimits_{k<n}a_{k}x_{j_{k}}\right\| .
\end{equation*}

\item  $C$-unconditional, where $C<\infty $ is a constant, if 
\begin{equation*}
C^{-1}\left\| \sum\nolimits_{k<n}a_{k}\epsilon _{k}x_{i_{k}}\right\| \leq
\left\| \sum\nolimits_{k<n}a_{k}x_{i_{k}}\right\| \leq C\left\|
\sum\nolimits_{k<n}a_{k}\epsilon _{k}x_{i_{k}}\right\|
\end{equation*}
for any choosing of $n<\omega $; $\{a_{k}:k<n\}$; $i_{0}<i_{1}<...<i_{n-1}<%
\alpha $ and of signs $\{\epsilon _{k}\in \{+,-\}:k<n\}$.

\item  Unconditional, if it is $C$-unconditional for some $C<\infty $.

\item  Symmetric, if for any $n<\omega $, any finite subset $I\subset 
\mathbb{\alpha }$ of cardinality $n$, any rearrangement $\varsigma $ of
elements of $I$ and any scalars $\{a_{i}:i\in I\}$, 
\begin{equation*}
\left\| \sum\nolimits_{i\in I}a_{i}z_{i}\right\| =\left\|
\sum\nolimits_{i\in I}a_{\varsigma (i)}z_{i}\right\| .
\end{equation*}

\item  Subsymmetric, if it is both spreading and 1-unconditional.
\end{itemize}
\end{definition}

Let $C<\infty $ be a constant. Two $\alpha $-sequences $\{x_{\beta }:\beta
<\alpha \}$ and $\{y_{\beta }:\beta <\alpha \}$ are said to be $C$\textit{%
-equivalent} if for any finite subset $I=\{i_{0}<i_{1}<...<i_{n-1}\}$ of $%
\alpha $ and for any choosing of scalars $\{a_{k}:k<n\}$ 
\begin{equation*}
C^{-1}\left\| \sum\nolimits_{k<n}a_{k}x_{i_{k}}\right\| \leq \left\|
\sum\nolimits_{k<n}a_{k}y_{i_{k}}\right\| \leq C\left\|
\sum\nolimits_{k<n}a_{k}x_{i_{k}}\right\| .
\end{equation*}

Two $\alpha $-sequences $\{x_{\beta }:\beta <\alpha \}$ and $\{y_{\beta
}:\beta <\alpha \}$ are said to be equivalent if they are $C$-equivalent for
some $C<\infty $.

\begin{remark}
The same definitions may be used in a case when instead of $\alpha $%
-sequences will be regarded families $\{x_{i}:i\in I\}\subset X$ indexed by
elements of a linearly ordered set $\left\langle I,\ll \right\rangle $. In a
such case it will be said about spreading families, unconditional families
and so on. $\omega $-sequences will be called sequences and may be denoted
like $\left( x_{n}\right) $.
\end{remark}

\section{Indices of divisibility}

\begin{definition}
Let $1\leq p\leq \infty $, $X^{f}$ be a class of finite equivalence. $X^{f}$
is said to be $p$-divisible if for some $Y\in X^{f}$ the space $l_{p}\left(
Y\right) $ belongs to $X^{f}$.
\end{definition}

\begin{remark}
It is easy to see that if $X^{f}$ is $p$-divisible then for every $Z\in
X^{f} $ a space $l_{p}\left( Z\right) $ belongs to $X^{f}$ too. Indeed,
since $l_{p}\left( Y\right) $ belongs to $X^{f}$, then for any ultrafilter $D
$ $\left( l_{p}\left( Y\right) \right) _{D}\in X^{f}$. Certainly, $%
l_{p}\left( \left( Y\right) _{D}\right) $ is isometric to a subspace of $%
\left( l_{p}\left( Y\right) \right) _{D}$. If $D$ is such that $Z$ is
isometric to a subspace of $\left( Y\right) _{D}$ then 
\begin{equation*}
l_{p}\left( Z\right) \hookrightarrow l_{p}\left( \left( Y\right) _{D}\right)
\hookrightarrow \left( l_{p}\left( Y\right) \right) _{D}
\end{equation*}
and, hence, $l_{p}\left( Z\right) \in X^{f}$.
\end{remark}

\begin{remark}
A simple criterion on $X^{f}$ to be $p$-divisible is: for any pair $A$, $%
B\in \frak{M}(X^{f})$ their $l_{p}$-sum $A\oplus _{p}B$ belongs to $\frak{M}%
(X^{f})$. Clearly, this condition satisfies when $X^{f}$ is $p$-divisible.
Conversely, if $A$, $B\in \frak{M}(X^{f})$ implies $A\oplus _{p}B\in \frak{M}%
(X^{f})$ for any $A$, $B$, then the space $W=l_{p}\left( A_{i},I\right) $,
where $\{A_{i}:i\in I\}$ is a numeration of all spaces from $\frak{M}(X^{f})$%
, belonging to $X^{f}$. Obviously, $l_{p}\left( W\right) \in X^{f}$ and,
thus, $X^{f}$ is $p$-divisible.
\end{remark}

Let $X$, $Y$ be Banach spaces. It will be said that $X$ is $Y$\textit{%
-saturated}, if any infinite-dimensional subspace of $X$ contains a
subspace, which is isomorphic to $Y$.

\begin{theorem}
Every $p$-divisible class $X^{f}$ contains a separable $l_{p}$-saturated
space.
\end{theorem}

\begin{proof}
Let $\left( A_{n}\right) _{n<\infty }$ be a dense subset of $\frak{M}(X^{f})$%
. Obviously, $l_{p}\left( A_{n}\right) \in X^{f}$ and is $l_{p}$-saturated.
\end{proof}

This simple result will play an important role in the following result.

\begin{definition}
Let $X^{f}$ be a class of finite equivalence. Its index of divisibility, $%
Index(X^{f})$ is a set of all such $p\in \lbrack 1.\infty ]$ that $X^{f}$ is 
$p$-divisible: 
\begin{equation*}
Index(X^{f})=\{p\in \lbrack 1.\infty ]:W\in X^{f}\Longrightarrow l_{p}\left(
W\right) \in X^{f}\}.
\end{equation*}
\end{definition}

For some classes $X^{f}$ their index of divisibility may be empty, $%
Index(X^{f})=\varnothing $.

E.g., $Index(\left( l_{p}\oplus _{2}l_{q}\right) ^{f})=\emptyset $ if $p\neq
q$.

Sometimes $Index(X^{f})$ consists of a single point: 
\begin{equation*}
Index(\left( l_{p}\right) ^{f})=\{p\}.
\end{equation*}

The maximal set $Index(X^{f})$ has the class $\left( l_{\infty }\right) ^{f}$%
: 
\begin{equation*}
Index(\left( l_{\infty }\right) ^{f})=[1,\infty ].
\end{equation*}

Indeed, the space $l_{p}\left( l_{\infty }\right) \in \left( l_{\infty
}\right) ^{f}$ for any $p$.

\begin{theorem}
If $Card\left( Index(X^{f})\right) \geq 2$ then the class $X^{f}$ contains a
space that has the Tsirelson property.
\end{theorem}

\begin{proof}
Let $p$, $q\in Index(X^{f})$; $p\neq q$. Then $X^{f}$ contains two separable
spaces: $X_{p}$, which is $p$-saturated and $X_{q}$, which is $q$-saturated.
Clearly, both and $X_{q}$ may be represented as closure of unions of chains
of their finite dimensional subspaces: 
\begin{eqnarray*}
X_{p} &=&\overline{\cup A_{k}};\text{ }A_{1}\hookrightarrow
A_{2}\hookrightarrow ...\hookrightarrow A_{n}\hookrightarrow ...; \\
X_{q} &=&\overline{\cup B_{k}};\text{ }B_{1}\hookrightarrow
B_{2}\hookrightarrow ...\hookrightarrow B_{n}\hookrightarrow ....
\end{eqnarray*}

Chose a sequence $\varepsilon _{n}\downarrow 0$ and define inductively a
sequence of isomorphic embeddings.

At the first step find the least number $n(1)$ such that there exists an
isomorphic embedding $u_{1}:A_{1}\rightarrow B_{n(1)}$ with $\left\|
u_{1}\right\| \left\| u_{1}^{-1}\right\| \leq 1+\varepsilon _{1}$; at the
same step choose the minimal $m(2)$ such that there exists an isomorphic
embedding $u_{2}:B_{n(1)}\rightarrow A_{m(2)}$ with $\left\| u_{2}\right\|
\left\| u_{2}^{-1}\right\| \leq 1+\varepsilon _{2}$.

If operators $u_{1}$, $u_{2}$, ..., $u_{2k}$ and numbers $n(1)$, $m(2)$, $%
n(3)$, ..., $m(2k)$ are already chosen, then $n(2k+1)$ is the least number
such that there exists an isomorphic embedding $u_{2k+1}:A_{m(2k)}%
\rightarrow B_{n(2k+1)}$ with $\left\| u_{2k+1}\right\| \left\|
u_{2k+1}^{-1}\right\| \leq 1+\varepsilon _{2k+1}$.

Also, chose $m(2k+2)$ as the least number such that there exists an
isomorphic embedding $u_{2k+2}:B_{n(2k+1)}\rightarrow A_{m(2k+2)}$ with $%
\left\| u_{2k+2}\right\| \left\| u_{2k+2}^{-1}\right\| \leq 1+\varepsilon
_{2k+2}$.

As a result it will be obtained a chain of finite dimensional spaces from $%
\frak{M}(X^{f})$: 
\begin{equation*}
cA_{1}\rightarrow B_{n(1)}\rightarrow A_{m(2)}\rightarrow ...\rightarrow
A_{m(2k)}\rightarrow B_{n(2k+1)}\rightarrow ...\text{,}
\end{equation*}
which may be regarded as a direct system.

Let $Y$ be a direct limit of a given direct system. Obviously, $Y$ may be
represented as the closure of a chain 
\begin{equation*}
A_{1}^{\prime }\hookrightarrow B_{n(1)}^{\prime }\hookrightarrow
A_{m(2)}^{\prime }\hookrightarrow ...\hookrightarrow A_{m(2k)}^{\prime
}\hookrightarrow B_{n(2k+1)}^{\prime }\hookrightarrow ...Y
\end{equation*}
of its finite dimensional subspaces, such that 
\begin{equation*}
d\left( A_{m(2k)}^{\prime },A_{m(2k)}\right) \leq 1+\varepsilon _{m\left(
2k\right) };\text{ \ }d\left( B_{n(2k+1)}^{\prime },B_{n(2k+1)}\right) \leq
1+\varepsilon _{m\left( 2k+1\right) }.
\end{equation*}

Assume that $l_{r}$ is isomorphic to a subspace of $Y$. Let $%
j:l_{r}\rightarrow Y$.

Then a chain 
\begin{equation*}
jl_{r}\cap A_{1}^{\prime }\hookrightarrow jl_{r}\cap B_{n(1)}^{\prime
}\hookrightarrow ...\hookrightarrow jl_{r}\cap A_{m(2k)}^{\prime
}\hookrightarrow jl_{r}\cap B_{n(2k+1)}^{\prime }\hookrightarrow ...\text{.}
\end{equation*}
contains a sub-chain 
\begin{equation*}
C_{1}\hookrightarrow C_{2}\hookrightarrow ...\hookrightarrow
C_{k}\hookrightarrow ...,
\end{equation*}
where $C_{k}$ is $\lambda $-isomorphic to a $l_{r}^{(s_{k})}$ for some $%
s_{k} $ and some $\lambda <\infty $, which does not depend on $k$.

However, this is impossible: if $\{jl_{r}\cap B_{n(2k+1)}^{\prime }:k<\infty
\}$ contains such a chain, then $l_{r}$ must be isomorphic to a subspace of $%
X_{q}$. If $\{jl_{r}\cap A_{m(2k)}^{\prime }:k<\infty \}$ contains such a
chain, then $l_{r}$ must be isomorphic to a subspace of $X_{p}$. Since $%
X_{p} $ is $p$-saturated and $X_{q}$ is $q$-saturated, this is impossible.
\end{proof}

\begin{corollary}
There exists a finite universal Banach space $Z$ that has the Tsirelson
property.
\end{corollary}

\begin{proof}
As was noted, $Index(\left( l_{\infty }\right) ^{f})=[1,\infty ]$ and,
hence, 
\begin{equation*}
Card\left( Index(\left( l_{\infty }\right) ^{f})\right) \geq 2.
\end{equation*}
\end{proof}

\begin{remark}
Recall that the original B. Tsirelson's space $[1]$ also was finitely
universal.
\end{remark}

\begin{remark}
It may be given a more short proof of the theorem 2.

Indeed, let $p$, $q\in Index(\left( X\right) ^{f})$. Let $\{A_{i}:i<\infty
\} $ be a dense in $\frak{M}(X^{f})$ sequence. Consider a sequence of spaces 
$\left( B_{k}\right) _{k<\infty }$ that are defined inductively: 
\begin{eqnarray*}
B_{1} &=&A_{1};\text{ }B_{2}=B_{1}\oplus _{p}A_{2};\text{ }B_{3}=B_{2}\oplus
_{q}A_{3};\text{ }... \\
B_{2k} &=&B_{2k-1}\oplus _{p}A_{2k};\text{ }B_{2k+1}=B_{2k}\oplus
_{q}A_{2k+1};\text{ }...
\end{eqnarray*}

It easy to see that $\underset{\rightarrow }{\lim }B_{n}$ has the Tsirelson
property.

Nevertheless, the proof of the theorem 2 (in difference from the given above
one) may be applied in more general cases, as it will be shown below.
\end{remark}

It may be presented other examples of classes $X^{f}$, such that $%
Index(X^{f})$ is more then one single point.

\begin{theorem}
For any closed subset $e\subset \lbrack 1,\infty )$ there exists a class $%
\left( X_{e}\right) ^{f}$ such that 
\begin{equation*}
Index(\left( X_{e}\right) ^{f})=e.
\end{equation*}
\end{theorem}

\begin{proof}
Let $\{p_{i}:i<\infty \}$ be a countable dense subset of $e$. Consider an
infinite matrix 
\begin{equation*}
\begin{array}{llll}
p_{1} & p_{2} & p_{3} & ... \\ 
p_{1} & p_{2} & p_{3} & ... \\ 
p_{1} & p_{2} & p_{3} & ... \\ 
... & ... & ... & ...
\end{array}
\end{equation*}

Let $\left( q_{i}\right) $ be an enumeration of it in a sequence, such that
each of $p_{k}$'s occurs among $q_{i}$'s infinitely many times. It may be
assumed that $p_{1}\neq 2$.

A class $\left( X_{e}\right) ^{f}$ will be constructed by induction. Let $%
X_{0}$ be a Banach space that generates a class $\left( X_{0}\right) ^{f}$
with $Index(\left( X_{0}\right) ^{f})=\varnothing $. Let $%
X_{0}=l_{p_{1}}\oplus _{2}l_{2}$ 
\begin{equation*}
X_{1}=l_{q_{1}}(X_{0});\text{ }X_{2}=l_{q_{2}}\left( X_{1}\right) ;\text{ }%
...;\text{ }X_{n+1}=l_{q_{n+1}}(X_{n});\text{ }...
\end{equation*}

Certainly, $X_{i}$ may be regarded as a subspace of $X_{i+1}$ and, thus, the
direct chain 
\begin{equation*}
X_{1}\hookrightarrow X_{2}\hookrightarrow \text{ }...\hookrightarrow
X_{n+1}\hookrightarrow ...
\end{equation*}
has a direct limit $X_{e}=\overline{\cup X_{i}}$. It is clear that the class 
$\left( X_{e}\right) ^{f}$ has desired properties.

Indeed, for any $A$, $B\in \frak{M}(\left( X_{e}\right) ^{f})$ and any $%
p_{i} $ their $l_{p_{i}}$-sum $A\oplus _{p_{i}}B$ belongs to $\frak{M}%
(\left( X_{e}\right) ^{f})$, as it follows from the construction. Hence, 
\begin{equation*}
\{p_{i}:i<\infty \}\subseteq Index(\left( X_{e}\right) ^{f}).
\end{equation*}

If $A\subseteq Index(\left( X_{e}\right) ^{f})$ then its closure $\overline{A%
}$ also belongs to $Index(\left( X_{e}\right) ^{f})$. This fact easily
follows from the closedness of $\frak{M}(\left( X_{e}\right) ^{f})$. If $%
p\notin e$ then $p\notin Index(\left( X_{e}\right) ^{f})$ since in a
contrary case any sum $A\oplus _{p}B$, where $A$, $B\in \frak{M}(\left(
X_{e}\right) ^{f})$ belongs to $\frak{M}(\left( X_{n}\right) ^{f})\ $for
some $n$. Certainly, this is impossible if $p\notin e$.
\end{proof}

\begin{remark}
Of course, $Index(X^{f})\subseteq S(X^{f})$ for any Banach space $X$. The
method, presented above allow to construct classes $X^{f}$ with $%
Index(X^{f})=S(X^{f})$.

E.g. for a two-point set $\{2,q\}$, where $q>2$ it may be considered:

as $X_{0}$ the space $l_{q}$ (for which $S(l_{q})=\{2,q\}$);

as $X_{1}$ - the space $l_{2}(l_{q})$;

$X_{3}=l_{q}\left( l_{2}(l_{q})\right) $;

$X_{4}=l_{2}(l_{q}(l_{2}(l_{q})))$ etc.

Clearly, the direct limit of $X_{i}$'s has the two-pointed $l_{p}$-spectrum
and the same index of divisibility.
\end{remark}

\section{Superstable classes of finite equivalence}

It was shown that some classes of finite equivalence contain spaces with the
Tsirelson properties.

From the other hand, there are classes $X^{f}$ that has ''anti-Tsirelson
property'': every representative of a such class contains some $l_{p}$. Some
of these classes may be pick out by using of \textit{stable Banach spaces},
which was introduced by J.-L. Krivine and B. Maurey [3] and their
generalization - \textit{superstable Banach spaces} that were defined by J.
Reynaud [4].

\begin{definition}
A Banach space $X$ is said to be stable provided for any two sequences $%
\left( x_{n}\right) $ and $\left( y_{m}\right) $ of its elements and every
pair of ultrafilter $D$, $E$ over $\mathbb{N}$ 
\begin{equation*}
\lim_{D\left( n\right) }\lim_{E\left( m\right) }\left\| x_{n}+y_{m}\right\|
=\lim_{E\left( m\right) }\lim_{D\left( n\right) }\left\| x_{n}+y_{m}\right\|
.
\end{equation*}
\end{definition}

The notations $D(n)$ and $E(m)$ are used here (instead of $D$ and $E$) to
underline the variable ($n$ or $m$ respectively) in expressions like $%
\lim_{D\left( n\right) }f(n,m)$.

\begin{definition}
(Cf. $[4]$). A Banach space $X$ is said to be superstable if every its
ultrapower $\left( X\right) _{D}$ is stable.
\end{definition}

\begin{theorem}
Let $X^{f}$ be a class of finite equivalence. $X^{f}$ contains a superstable
Banach space if and only if every space $Y\in X^{f}$ (and, as a consequence,
every space $W$, which is finitely representable in $X$) is stable.
\end{theorem}

\begin{proof}
If $Y\in X^{f}$ is superstable then each its subspace is stable because of
the property of a Banach space to be stable is inherited by its subspaces.
Hence, each subspace of every ultrapower $\left( Y\right) _{D}$ is stable
too, because of $\{Z:Z<_{f}Y\}$ is coincide with the set 
\begin{equation*}
\{Z:Z\text{ \ is isometric to a subspace of some ultrapower }\left( Y\right)
_{D}
\end{equation*}

Conversely, if every $Y\in X^{f}$ is stable, then all ultrapowers of $Y$ are
stable too.
\end{proof}

\begin{definition}
A class $X^{f}$ of finite equivalence that contains a superstable space will
be called a superstable class.
\end{definition}

In [3] it was shown:

\textit{Any stable Banach space }$X$\textit{\ is weakly sequentially
complete. }

\textit{Every subspace of }$X$\textit{\ contains a subspace isomorphic to
some} $l_{p}$ ($1\leq p<\infty $).

\begin{theorem}
Let $X^{f}$ be a superstable class. Then there exists such $p\in \lbrack
1,\infty )$ that every $Y\in X^{f}$ contains a subspace, which is isomorphic
to $l_{p}$.
\end{theorem}

\begin{proof}
Let $X\in \mathcal{B}$. Let 
\begin{equation*}
T(X)=\{p:l_{p}\text{ is isomorphic to a subspace of }X\}.
\end{equation*}
Let $Y$, $Z\in X^{f}$. Let both $Y$ and $Z$ are separable and stable. Thus, $%
T(Y)\neq \varnothing $ and $T(Z)\neq \varnothing $. Assume that $T(Y)\cap
T(Z)=\varnothing $. Similarly to the proof of the theorem 2 may be
constructed a space $W\in X^{f}$ that has the Tsirelson property. However
this conflicts with the superstability of $X^{f}$. Hence $\cap \{T(Y):Y\in
X^{f}\}\neq \varnothing $.
\end{proof}

\begin{remark}
There are classes $X^{f}$ such that the intersection $\cap \{T(Y):Y\in
X^{f}\}$ consists exactly of $n$ points.

Indeed, let $X=l_{q_{1}}\oplus l_{q_{2}}\oplus ...\oplus l_{q_{n}}$, where $%
2<q_{1}<q_{2}<...<q_{n}$. Clearly, 
\begin{equation*}
\cap \{T(Y):Y\in X^{f}\}=\{q_{1},q_{2},...,q_{n}\}.
\end{equation*}
\end{remark}

\begin{problem}
Whether there exists such a class $X^{f}$ that $\cap \{T(Y):Y\in X^{f}\}$ is
infinite?
\end{problem}

In what follows will be needed a definition.

\begin{definition}
(Cf. $[5]$) Let $X$ be a Banach space, $\left( x_{n}\right) \subset X$ be a
nontrivial normed sequence of elements of $X$ (i.e. $\left( x_{n}\right) $
contains no Cauchy subsequences); $D$ be an ultrafilter over $\mathbb{N}$.
For a finite sequence $\left( a_{i}\right) _{i=1}^{n}\subset \mathbb{R}^{n}$
let 
\begin{equation*}
l\left( \left( a_{i}\right) _{i=1}^{n}\right) =\lim_{D\left( m_{1}\right)
}\lim_{D\left( m_{2}\right) }...\lim_{D\left( m_{n}\right) }\{\left\|
\sum\nolimits_{k=1}^{n}a_{k}x_{m_{k}}\right\| :m_{1}<m_{2}<...<m_{n}\}.
\end{equation*}

Let $sm(X,\left( x_{n}\right) ,D)$ be a completition of a linear space $%
c_{00}$ of all sequences $\left( a_{i}\right) _{i=1}^{\infty }$ of real
numbers such that all but finitely many $a_{i}$'s are equal to zero.

The space $sm(X,\left( x_{n}\right) ,D)$ is said to be a spreading model of
the space $X$, which is based on the sequence $\left( x_{n}\right) $ and on
the ultrafilter $D$.
\end{definition}

Clearly, any spreading model of a given Banach space $X$ has a spreading
basis and is finitely representable in $X$.

In [3] it was shown that every spreading model of a stable Banach space $X$
has a symmetric basis. The following result shows that in a superstable
class X$^{f}$ the converse is also true.

\begin{definition}
Let $X$ be a Banach space. Its $IS$-spectrum $IS(X)$ is a set of all
(separable) spaces $\left\langle Y,\left( y_{i}\right) \right\rangle $ with
a spreading basis $\left( y_{i}\right) $ which are finitely representable in 
$X$.
\end{definition}

Notice that if $\left( y_{i}^{\prime }\right) $ and $\left( y_{i}^{\prime
\prime }\right) $ are different spreading bases of $Y$ then $\left\langle
Y,\left( y_{i}^{\prime }\right) \right\rangle $ and $\left\langle Y,\left(
y_{i}^{\prime \prime }\right) \right\rangle $ are regarded as different
members of $IS(X)$.

The last reservation may be omitted if one assumes that members of $IS(X)$
are nonseparable spaces $\left\langle Y,\left( y_{\alpha }\right) _{\alpha
<\omega _{1}}\right\rangle $ with uncountable spreading bases $\left(
y_{\alpha }\right) _{\alpha <\omega _{1}}$ (here and below $\omega _{1}$
denotes the first uncountable cardinal). The proof of this assertion will be
given below.

\begin{theorem}
A class $X^{f}$ is superstable if and only if every member $\left\langle
Y,\left( y_{i}\right) \right\rangle $ of its $IS$-spectrum has a symmetric
basis.
\end{theorem}

\begin{proof}
Let $X^{f}$ be superstable; $Y\in X^{f}$. According to [3] every spreading
model of $Y$ has a symmetric basis. Let $\left( Y\right) _{D}$ be an
ultrapower by a countably incomplete ultrafilter. Then (see [6]) $\left(
Y\right) _{D}$ contains any separable Banach space which is finitely
representable in $Y$. In particular, any space $\left\langle Z,\left(
z_{i}\right) \right\rangle $ of $IS(X^{f})$ is isometric to a subspace of $%
\left( Y\right) _{D}$. Since $\left( z_{i}\right) $ is a spreading sequence, 
$sm(\left( Y\right) _{D},\left( z_{i}\right) ,E)$ is isometric to $%
\left\langle Z,\left( z_{i}\right) \right\rangle $ for any ultrafilter $E$.
Hence, $\left( z_{i}\right) $ is a symmetric sequence.

Conversely, assume that every $\left\langle Z,\left( z_{i}\right)
\right\rangle $ has a symmetric basis. Assume that $X$ is not superstable.
Then there exists a space from $X^{f}$ which is not stable (it may be
assumed that $X$ is not stable itself). According to [3] there are such
sequences $\left( x_{n}\right) $ and $\left( y_{m}\right) $ of elements of $%
X $ that 
\begin{equation*}
\sup_{m<n}\left\| x_{n}+y_{m}\right\| >\inf_{m>n}\left\| x_{n}+y_{m}\right\|
.
\end{equation*}

Let $D$ be a countably incomplete ultrafilter over $\mathbb{N}$. Let $X_{0}%
\overset{def}{=}X$; 
\begin{equation*}
X_{n}\overset{def}{=}\left( X_{n-1}\right) _{D};\text{ }n=1,\text{ }2,\text{ 
}...;\text{ \ }X_{\infty }=\overline{\cup _{n\geq 1}X_{n}}.
\end{equation*}
Here is assumed that $X_{n}$ is a subspace of $X_{n+1}=\left( X_{n}\right)
_{D}$ under the canonical embedding $d_{X_{n}}:X_{n}\rightarrow \left(
X_{n}\right) _{D}$.

Let $D$, $E$ be ultrafilters over $\mathbb{N}$. Their product $D\times E$ is
a set of all subsets $A$ of $\mathbb{N}\times \mathbb{N}$ that are given by 
\begin{equation*}
\{j\in \mathbb{N}:\{i\in \mathbb{N}:\left( i,j\right) \in A\}\in D\}\in E.
\end{equation*}
Certainly, $D\times E$ is an ultrafilter and for every Banach space $Z$ the
ultrapower $\left( Z\right) _{D\times E}$ may be in a natural way identified
with $\left( \left( Z\right) _{D}\right) _{E}$.

So, the sequence $\left( x_{n}\right) \subset X$ defines elements 
\begin{eqnarray*}
\frak{x}_{1} &=&\left( x_{n}\right) _{D}\in \left( X\right) _{D}; \\
\frak{x}_{2} &=&\left( x_{n}\right) _{D\times D}\in \left( \left( X\right)
_{D}\right) _{D}; \\
&&...\text{ \ \ \ \ }...\text{ \ \ \ \ }...\text{ \ \ \ \ }...\text{ \ \ \ \ 
}...\text{ \ \ \ \ }... \\
\frak{x}_{k} &=&\left( x_{n}\right) \underset{k\text{ }times}{_{\underbrace{%
D\times D\times ...\times D}}}\in \underset{k\text{ }times}{\underbrace{%
\left( \left( \left( X\right) _{D}\right) _{D}...\right) _{D}};} \\
&&...\text{ \ \ \ \ }...\text{ \ \ \ \ }...\text{ \ \ \ \ }...\text{ \ \ \ \ 
}...\text{ \ \ \ \ }...
\end{eqnarray*}

Notice that $\frak{x}_{k}\in X_{k}\backslash X_{k-1}$. It is easy to verify
that $\left( \frak{x}_{k}\right) _{k<\infty }\subset X_{\infty }$ is a
spreading sequence. Since $X_{\infty }\in X^{f}$, it is symmetric. Moreover,
for any $z\in X$, where $X$ is regarded as a subspace of $X_{\infty }$ under
the direct limit of compositions 
\begin{equation*}
d_{X_{n}}\circ d_{X_{n-1}}\circ ...\circ d_{X_{0}}:X\rightarrow X_{n},
\end{equation*}
the following equality is satisfied: for any pair $m$, $n\in \mathbb{N}$ 
\begin{equation*}
\left\| x_{n}+z\right\| =\left\| x_{m}+z\right\| .
\end{equation*}
Since $\left( x_{n}\right) $ and $z$ are arbitrary elements of $X$, this
contradicts to the inequality $\sup_{m<n}\left\| x_{n}+y_{m}\right\|
>\inf_{m>n}\left\| x_{n}+y_{m}\right\| $.
\end{proof}

This result may be generalized to classes of crudely finite equivalence.

\begin{definition}
Let $X$, $Y$ be Banach spaces. $X$ is said to be crudely finitely
representable in $Y$ (in symbols: $X<_{F}Y$) if $X$ is isomorphic to a some
space which is finite representable in $Y$ ( it is easy to see that this
definition is equivalent to the standard one)
\end{definition}

Let $X\approx Y$ denotes that Banach spaces $X$ and $Y$ are isomorphic.

A class $X^{F}$ of crudely finite equivalence, which is generated by a
Banach space $X$ is given by 
\begin{equation*}
X^{F}=\{Y\in \mathcal{B}:Y<_{F}X\text{ \ and }X<_{F}Y\}=\cup
\{Y^{f}:Y\approx X\}.
\end{equation*}

\begin{definition}
A class $X^{F}$ of crudely finite equivalence is said to be crudely
superstable if it contains a superstable space.
\end{definition}

In other words, $X^{F}$ is crudely superstable if one of classes $Y^{f}$
that form $X^{F}$ (i.e one of classes of the union $\cup \{Y^{f}:Y\approx
X\}=X^{F}$) is superstable.

Certainly, any crudely superstable class $X^{F}$ has the property:

\textit{Every Banach space }$Y$\textit{, which is finitely representable in
a some space }$Z\in X^{F}$\textit{\ contains a subspace that is isomorphic
to some }$l_{p}$\textit{\ (}$1<p<\infty $\textit{)}.

\begin{theorem}
A class $X^{F}$ is crudely superstable if and only if for every $Z\in X^{F}$
its $IS$-spectrum $IS(Z)$ consists of spaces $\left\langle W,\left(
w_{n}\right) \right\rangle $, whose natural bases $\left( w_{n}\right) $ are 
$c_{Z}$-equivalent to symmetric bases where the constant $c_{Z}$ depends
only on $Z$.
\end{theorem}

\begin{proof}
Let $X^{F}$ be crudely superstable. Then some $Y\in X^{F}$ is superstable
and for any $Z$, which is crudely finitely representable in $Y$, and for any
space $\left\langle W,\left( w_{n}\right) \right\rangle $ from the $IS$%
-spectrum $IS(Z)$, its natural basis $\left( w_{n}\right) $ is $d(Z,Y_{1})$%
-equivalent to a symmetric one, where $Y_{1}\in Y^{f}$.

Conversely, let $X_{0}$ be such that every space $\left\langle W,\left(
w_{n}\right) \right\rangle \in IS(X_{0})$ has a basis $\left( w_{n}\right) $%
, which is equivalent to a symmetric one (certainly, this is equivalent to
the assertion that for every $Z\in X^{F}$ its $IS$-spectrum $IS(Z)$ consists
of spaces $\left\langle W^{\prime },\left( w_{n}^{\prime }\right)
\right\rangle $, whose natural bases $\left( w_{n}^{\prime }\right) $ are
equivalent to symmetric bases. It is easy to show that there exists a
constant $c_{X}$ such that every space $\left\langle W,\left( w_{n}\right)
\right\rangle \in IS(X_{0})$ has a basis $\left( w_{n}\right) $ which is $%
c_{X}$- equivalent to a symmetric one.

Indeed, let $\left( c_{k}\right) $ be a sequence of real numbers with a
property: for every $k<\infty $ there exists $\left\langle W_{k},\left(
w_{n}^{k}\right) \right\rangle \in IS(X_{0})$ such that $\left(
w_{n}^{k}\right) $ is $c_{k}$-equivalent to a symmetric basis and is not $%
c_{k-1}$-equivalent to any symmetric basis. Without loss of generality it
may be assumed that all spaces $\left( W_{k}\right) $ are subspaces of a
space from $\left( X_{0}\right) ^{f}$, e.g. $W_{k}\hookrightarrow X_{0}$.
Consider an ultrapower $\left( X_{0}\right) _{D}$ and its elements $\frak{w}%
_{k}=\left( w_{n}^{k}\right) _{D(k)}$. Clearly $\left( \frak{w}_{k}\right)
_{k<\infty }\subset \left( X_{0}\right) _{D}$ is a spreading sequence that
is not equivalent to any symmetric sequence.

Consider some $Z\in X^{F}$, such that $Z$ contains any space from its $IS$%
-spectrum. Let $c_{Z}$ be the corresponding constant which was defined
above. Let $\{Z_{\alpha }:\alpha <\varkappa \}$ ($\varkappa $ is a cardinal
number) be a numeration of all subspaces of $Z$ that may be represented as $%
span(z_{i}^{(\alpha )})$ (i.e. as a closure of linear span of $%
\{z_{i}^{(\alpha )}:i<\infty \}$) for such sequences $\{z_{i}^{(\alpha
)}:i<\infty \}$ (not necessary spreading ones) that are $c_{Z}$-equivalent
to symmetric sequences.

Using the standard procedure of renorming, due to A. Pe\l czy\'{n}ski [12],
it may be constructed a space $Z_{\infty }\approx Z$ such that $Z_{\infty }$
contains as a subspace every space $W_{\infty }$ from $IS(Z_{\infty })$,
which (by the renorming procedure) has a symmetric basis.

From the theorem 6 follows that $Z_{\infty }$ (and, hence, the whole class $%
\left( Z_{\infty }\right) ^{f}$) is superstable. Since $Z_{\infty }\in X^{F}$
the class $X^{F}$ is crudely superstable.

Indeed, it is sufficient to choose as a unit ball $B(Z_{\infty })=\{w\in
Z_{\infty }:\left\| w\right\| \leq 1\}$ a convex hull of the union of a set $%
B(Z)$ with sets $\{c_{Z}^{-1}j_{\alpha }B(i_{\alpha }Z_{\alpha }):$ $\alpha
<\varkappa \}$, where $i_{\alpha }:Z_{\alpha }\rightarrow W_{\alpha }$ is an
isomorphism between $Z_{\alpha }$ and a space $W_{\alpha }$ with a symmetric
basis $\left( w_{n}^{\left( \alpha \right) }\right) $, which is given by $%
i_{\alpha }\left( z_{n}^{(\alpha )}\right) =w_{n}^{\left( \alpha \right) }$
for all $n<\infty $; $\left\| i_{\alpha }\right\| \left\| i_{\alpha
}^{-1}\right\| \leq c_{Z}$; $j_{\alpha }$ is an embedding of the unit ball $%
B(W_{\alpha })$ in a set $c_{Z}B(Z_{\alpha })$, which is given by $%
c_{Z}B(Z_{\alpha })=\{z\in Z_{\alpha }:c_{Z}^{-1}z\in B\left( Z_{\alpha
}\right) \}$. Namely, 
\begin{equation*}
B(Z_{\infty })=conv\{B(Z)\cup \left( \cup \{c_{Z}^{-1}j_{\alpha }B(i_{\alpha
}Z_{\alpha }):\alpha <\varkappa \}\right) \}.
\end{equation*}
\end{proof}

\section{Spaces spanned by spreading uncountable sequences}

From results of the previous section follows that spreading sequences play
an important role in the study of the Tsirelson property.

Ordinals will be understands in the von Neuman sense: an ordinal $\alpha $
will be regarded as a set $\{\beta :\beta <\alpha \}$ of all ordinals that
are less then $\alpha $.

Small Greece letters $\alpha $, $\beta $, $\gamma $, $\delta $, $\zeta $
denote \textit{ordinals}; $\varkappa $, $\tau $ be \textit{cardinals}. In
what follows cardinals will be identified with initial ordinals. Finite
ordinals and cardinals may be denoted also by small Latin letters $i$, $j$, $%
k$, $m$, $n$.

The least infinite ordinal (cardinal) is denoted by $\omega $; the first
uncountable ordinal (cardinal) - by $\omega _{1}$.

For $\alpha $, $\beta $ are\ ordinals, the symbol $\alpha ^{\beta }$ denotes
the \textit{ordinal degree}.

A symbol $2^{\tau }$ denotes the \textit{cardinal degree} - a cardinality of
a set $\exp \tau $\ of all subsets of a cardinal $\tau $.

Let $X$ be a Banach space, $A\subset X$.

The \textit{linear span} $lin(A)$, i.e. a set of all linear combinations of
elements of $A$ need not to be a closed subspace of $X$. A closure $%
\overline{lin(A)}$ will be denoted by $span(A)$.

A \textit{dimension} of a Banach space $X$, $dim(X)$ is the least
cardinality of a subset $A\subset X$ such that $span(A)=X$.

Let $\{x_{n}:n<\omega \}$ be a spreading sequence.

It generates a sequence $\mathcal{N}=\left( \mathcal{N}_{m}\right)
_{m<\omega }$ of norms: $\mathcal{N}_{m}$ is a norm on $\mathbb{R}^{m}$,
which is given by 
\begin{equation*}
\mathcal{N}_{m}\left( \left( a_{i}\right) _{i=1}^{m}\right) =\left\|
\sum\nolimits_{i=1}^{m}a_{i}x_{i}\right\| .
\end{equation*}

Let $\left\langle I,\ll \right\rangle $ be a linearly ordered set.

Consider a vector space $c_{00}\left( I,\ll \right) $ of all families $%
\{a_{\iota }:\iota \in I\}$ of scalars all but finitely many elements of
which are vanished.

$c_{00}\left( I,\ll \right) $ may be equipped with a norm $\mathcal{N}$
which is given by 
\begin{equation*}
\mathcal{N}\left( \{a_{\iota }:\iota \in I\}\right) =\mathcal{N}_{m}\left(
a_{i_{1}}\text{, }a_{i_{2}}\text{, ...,}a_{i_{m}}\right) ,
\end{equation*}
where $m=card\{\iota \in I:a_{\iota }\neq 0\}$; $i_{1}\ll i_{2}\ll ...\ll
i_{m}$ and all $\{a_{i_{k}}:k=1$, $2$, ...,$m\}$ are differ from $0$.

The completition of the normed space $\left\langle c_{00}\left( I,\ll
\right) ,\mathcal{N}\right\rangle $ is a Banach space, which will be denoted
by $X_{\mathcal{N}}(I,\ll )$.

So, the spreading sequence $\{x_{n}:n<\omega \}=\frak{x}$ generates a class
of Banach spaces of kind $X_{\mathcal{N}}(I,\ll )$, which will be called 
\textit{the tower }$\left\lfloor \frak{x}\right\rfloor $\textit{, generated
by }$\frak{x}$\textit{.}

The following result is obvious ant its proof is omitted.

\begin{theorem}
Let $\{x_{n}:n<\omega \}=\frak{x}$ be a spreading sequence; $\left\lfloor 
\frak{x}\right\rfloor $ be the corresponding tower.

If $\frak{x}$ is symmetric, then all spaces $Y\in \left\lfloor \frak{x}%
\right\rfloor $ of the dimension $\tau $ are pairwice isometric, where $\tau 
$ is an arbitrary infinite cardinal.

If $\frak{x}$ is equivalent to a symmetric sequence then all $Y\in
\left\lfloor \frak{x}\right\rfloor $ of the dimension $\tau $ are pairwice
isomorphic.
\end{theorem}

Let 
\begin{equation*}
\mathcal{B}_{\varkappa }=\{X\in \mathcal{B}:dim(X)=\varkappa \};
\end{equation*}
\begin{equation*}
X^{\approx }=\{Y\in \mathcal{B}:Y\approx X\}.
\end{equation*}

Let $\mathcal{K}$ be a class of Banach spaces, which is closed under
isomorphisms (i.e. $X\in \mathcal{K}$ implies that $X^{\approx }\subset 
\mathcal{K}$. Let 
\begin{equation*}
\mathcal{K}^{\approx }=\{X^{\approx }:X\in \mathcal{K.}
\end{equation*}

It may be shown that if $\frak{x}$ is not symmetric (resp., is not
equivalent to a symmetric sequence) then the cardinality $card(\left\lfloor 
\frak{x}\right\rfloor \cap \mathcal{B}_{\varkappa })=2^{\varkappa }$ (resp.,
the cardinality $card\left( (\left\lfloor \frak{x}\right\rfloor \cap 
\mathcal{B}_{\varkappa })^{\approx }\right) =2^{\varkappa }$).

Let $\varkappa $ be a cardinal; $\sigma :\varkappa \rightarrow \varkappa $
be a transposition (i.e. one-to-one mapping of $\varkappa $ onto $\varkappa $%
).

It will be said that $\sigma $ is almost identical if the correlation 
\begin{equation*}
\gamma _{1}<\gamma _{2}\Rightarrow \sigma \gamma _{1}<\sigma \gamma _{2}%
\text{, where }\gamma _{1}\text{, }\gamma _{2}<\varkappa
\end{equation*}
may get broken at most finitely many times.

\begin{lemma}
Let $\{x_{n}:n<\omega \}$ be a spreading sequence that is not symmetric.

Let $\sigma :\omega \rightarrow \omega $ be a transposition, which is not
almost identical.

If sequences $\{x_{n}:n<\omega \}$ and $\{x_{\sigma n}:n<\omega \}$ are
equivalent then both of them are equivalent to a symmetric sequences.
\end{lemma}

\begin{proof}
Consider a sequence $\{y_{\alpha }:\alpha <\omega ^{2}\}$, which is given by 
\begin{equation*}
\{y_{\omega \cdot 2k+n}=x_{n};\text{ \ }y_{\omega \cdot \left( 2k+1\right)
+n}=x_{\sigma n}:k<\omega ;\text{ }n<\omega \}.
\end{equation*}
This sequence is equivalent to a sequence $\{y_{\alpha }^{\prime }:\alpha
<\omega ^{2}\}$ that belongs to the tower $\left\lfloor \frak{x}%
\right\rfloor $, generated by $\{x_{n}:n<\omega \}=\frak{x}$.

Certainly, this is equivalent to 
\begin{equation*}
C^{-1}\left\| \sum\nolimits_{k<n}a_{k}x_{i_{k}}\right\| \leq \left\|
\sum\nolimits_{k<n}a_{k}y_{\alpha _{k}}\right\| \leq C\left\|
\sum\nolimits_{k<n}a_{k}x_{i_{k}}\right\|
\end{equation*}
for every $n<\omega $; every scalars $\left( a_{k}\right) _{k<n}$ every
choice $\alpha _{0}<\alpha _{1}<...<\alpha _{n-1}<\omega ^{2}$ and $%
i_{0}<i_{1}<...<i_{n-1}<\omega $ and some $C<\infty $ that depends only on
equivalence constant between $\{x_{n}:n<\omega \}$ and $\{x_{\sigma
n}:n<\omega \}$. Since $\sigma $ has only a finite number of inversions, our
definition of $\{y_{\alpha }:\alpha <\omega ^{2}\}$ implies that $\frak{x}$
is equivalent to a symmetric sequence
\end{proof}

\begin{theorem}
Let $\alpha $, $\beta $ be ordinals; $\omega \leq \beta \leq \alpha $. Let $%
\{x_{\gamma }:\gamma <\alpha \}$ be a subsymmetric $\alpha $-sequence which
is not equivalent to any $\alpha $-symmetric sequence. Let $X_{\alpha
}=span(\{x_{\gamma }:\gamma <\alpha \})$; $X_{\beta }=span(\{x_{\gamma
}:\gamma <\beta \})$. If $\beta ^{\omega }<\alpha $ then spaces $X_{\alpha }$
and $X_{\beta }$ are not isomorphic.
\end{theorem}

\begin{proof}
Assume that $X_{\alpha }$ is isomorphic to a subspace $Z$\ of $X_{\beta }$.

Let $I:X_{\alpha }\rightarrow X_{\beta }$ be the corresponding operator of
isomorphic embedding. Without loss of generality it may be assumed that an
image of element $x_{\gamma }$ ($\gamma <\alpha $) in $X_{\beta }$ is a
finite linear combination with rational coefficients of some $x_{\zeta }$'s (%
$\zeta <\beta $): 
\begin{equation*}
Ix_{\gamma }=\sum\nolimits_{k=0}^{n(\gamma )}a_{k}^{\gamma }x_{\zeta
_{k}\left( \gamma \right) };\text{ \ }\zeta _{0}\left( \gamma \right) <\zeta
_{1}\left( \gamma \right) <...<\zeta _{n(\gamma )}\left( \gamma \right)
<\beta .
\end{equation*}
Thus, to any $x_{\gamma }$ corresponds a finite sequence of rational numbers 
$\left( a_{k}^{\gamma }\right) _{k=0}^{n(\gamma )}$.

Let $\left( p_{n}\right) _{n<\omega }$ be a numeration of all finite
sequences of rationals.

The set $\left( p_{n}\right) _{n<\omega }$ generates a partition of $\alpha $
on parts $\left( P_{n}\right) _{n<\omega }$ in a following way: an ordinal $%
\gamma <\alpha $ belongs to $P_{n}$ if $\left( a_{k}^{\gamma }\right)
_{k=0}^{n(\gamma )}=p_{n}$.

Since $\alpha >\beta ^{\omega }\geq \omega ^{\omega }$, one of $P_{n}$'s
should contain a sequence $\{\gamma _{i}:i<\delta \}$ of ordinals, which
order type $\delta $ (in a natural order : $i<j$ implies that $\gamma
_{i}<\gamma _{j}$) is greater then $\beta ^{\omega }$. Clearly, for all such 
$\gamma _{i}$, 
\begin{equation*}
Ix_{\gamma _{i}}=\sum\nolimits_{k=0}^{n}a_{k}x_{\zeta _{k}\left( \gamma
_{i}\right) },
\end{equation*}
where $n$ and $\left( a_{k}\right) _{k=1}^{n}$ do not depend on $i$.

A set of all sequences $\zeta _{0}\left( \gamma _{i}\right) <\zeta
_{1}\left( \gamma _{i}\right) <...<\zeta _{n}\left( \gamma _{i}\right)
<\beta $ can be ordered to have only the order type $\leq \beta ^{n}$.
Hence, conditions $\gamma _{i_{1}}<\gamma _{i_{2}}\Rightarrow \zeta
_{0}\left( \gamma _{i_{1}}\right) <\zeta _{0}\left( \gamma _{i_{2}}\right) $
must get broken for infinite many pairs $\gamma _{i_{1}}$, $\gamma _{i_{2}}$.

The inequality 
\begin{equation*}
C^{-1}\left\| \sum\nolimits_{k<m}b_{k}x_{\gamma _{k}}\right\| \leq \left\|
\sum\nolimits_{k<m}b_{k}\left( \sum\nolimits_{k=0}^{n}a_{k}x_{\zeta
_{k}\left( \gamma _{i}\right) }\right) \right\| \leq C\left\|
\sum\nolimits_{k<m}b_{k}x_{\gamma _{k}}\right\| ,
\end{equation*}
where $C=d(X_{\alpha },Z)$, shows that $\{x_{\gamma }:\gamma <\alpha \}$ is
equivalent to a symmetric sequence.
\end{proof}

\begin{theorem}
Let $\frak{x}=\{x_{n}:n<\omega \}$ be a subsymmetric sequence, which is not
isomorphic to a symmetric one. Then for any cardinal $\varkappa \geq \omega $
there exists $2^{\varkappa }$ pairwice non isomorphic Banach spaces of
dimension $\varkappa $, which belongs to the same tower $\left\lfloor \frak{x%
}\right\rfloor $.
\end{theorem}

\begin{proof}
Let $\{x_{\gamma }:\gamma <\varkappa \}$ be a subsymmetric $\varkappa $%
-sequence. Let $I=\left\langle I,\ll \right\rangle $ be a linearly ordered
set of cardinality $\varkappa $; $J=\left\langle I,<^{\prime }\right\rangle $
- another linear ordering of $I$.

Consider families $\{x_{i}:i\in I\}$ and $\{x_{j}:j\in J\}$ that are indexed
(and ordered) by elements of $I$ and $J$ respectively.

Let $X_{I}=span\{x_{i}:i\in I\}$; $X_{J}=span\{x_{j}:j\in J\}$.Certainly, $%
X_{I}$ and $X_{J}$ are isomorphic only if there exists one-to-one mappings
of embedding $u:I\rightarrow J$ and $w:J\rightarrow I$, which are almost
monotone in a following sense: 
\begin{eqnarray*}
i_{1} &\ll &i_{2}\Rightarrow u\left( i_{1}\right) <^{\prime }u\left(
i_{i}\right) \text{ \ for all but finitely many pairs }i_{1},i_{2}\in I; \\
j_{1} &<&^{\prime }j_{2}\Rightarrow w\left( j_{1}\right) \ll w\left(
j_{i}\right) \text{ \ for all but finitely many pairs }j_{1},j_{2}\in J.
\end{eqnarray*}

Since there exists $2^{\varkappa }$ orderings of $I$ for any pair of which
such almost monotone mappings do not exist, this prove the theorem.
\end{proof}

\begin{corollary}
The cardinality of the set of all Banach spaces of dimension $\varkappa $ is 
\begin{equation*}
card\left( \mathcal{B}_{\varkappa }\right) =2^{\varkappa }.
\end{equation*}
\end{corollary}

\begin{proof}
The inequality $card\left( \mathcal{B}_{\varkappa }\right) \geq 2^{\varkappa
}$ follows from the previous theorem. The inverse inequality is obvious.
\end{proof}

\begin{remark}
It is of interest that different sets $\mathcal{B}_{\varkappa }$ and $%
\mathcal{B}_{\tau }$ may be of the same cardinality. The appearance of a
such case depends on the model of the set theory that lies in the base of
all functional analysis.

E.g., if one assume the Martin axiom $MA$ with the negation of continuum
hypothesis $\urcorner CH$ then for all cardinals $\varkappa $ such that $%
\omega <\varkappa <2^{\omega }$ 
\begin{equation*}
card\mathcal{B}_{\varkappa }=card\mathcal{B}_{\omega }=2^{\omega }.
\end{equation*}
\end{remark}

It may be given an interesting cardinal criterion of superstability.

\begin{theorem}
Let $X$ be a Banach space; $X^{f}$ and $X^{F}$ be corresponding classes of
finite and crudely finite equivalence respectively, which are generated by $%
X $.

If $card\left( X^{f}\cap \mathcal{B}_{\varkappa }\right) <2^{\varkappa }$
then $X^{f}$ is a superstable class.

If $card\left( \left( X^{F}\cap \mathcal{B}_{\varkappa }\right) ^{\approx
}\right) <2^{\varkappa }$ then $X^{F}$ is crudely superstable.
\end{theorem}

\begin{proof}
According to [8] if $IS(X)$ contains a space $\left\langle W,\left(
w_{n}\right) \right\rangle $ with a spreading basis then for every cardinal $%
\varkappa $ and every transposition $\sigma $ of $\tau $ there are spaces $X$
and $X_{\sigma }$, which are finitely equivalent to $X$ and such that $%
dim(X)=dim(X_{\sigma })=\varkappa $; $X$ contains a subspace isometric to $%
W_{\varkappa }=span\{w_{\alpha }:\alpha <\varkappa \}$; $X_{\sigma }$
contains a subspace isometric to $W_{\sigma \varkappa }=span\{w_{\alpha
}:\alpha <\varkappa \}$. Improving arguments of [8] it may be shown that $X$
and $X_{\sigma }$ are isometric (resp., isomorphic) if and only if $%
W_{\varkappa }$ and $W_{\sigma \varkappa }$ are isometric (resp.,
isomorphic). Hence, if $X$ is not superstable then $card\left( X^{f}\cap 
\mathcal{B}_{\varkappa }\right) =2^{\varkappa }$. Similarly for the second
part of the theorem.
\end{proof}

It is known that a separable Banach space may have a lot of pairwice non
equivalent symmetric bases (e.g. the classical A. Pe\l czy\'{n}ski's space $%
P_{U}$, complementably universal in the class of all Banach spaces with
unconditional bases [9] has a continuum number of pairwice non equivalent
symmetric bases; cf. [10]).

However, if a nonseparable Banach space has a symmetric or subsymmetric
(uncountable) basis, this basis is unique up to equivalence.

\begin{theorem}
Let $\varkappa >\omega $ be a cardinal; $\{x_{\alpha }:\alpha <\varkappa \}$
and $\{y_{\alpha }:\alpha <\varkappa \}$ be subsymmetric sequences; $%
X=span\{x_{\alpha }:\alpha <\varkappa \}$; $Y=span\{y_{\alpha }:\alpha
<\varkappa \}$. If spaces $X$ and $Y$ are isomorphic then $\varkappa $%
-sequences $\{x_{\alpha }:\alpha <\varkappa \}$ and $\{y_{\alpha }:\alpha
<\varkappa \}$ are equivalent.
\end{theorem}

\begin{proof}
Let $u:X\rightarrow Y$ be an isomorphism; $\left\| u\right\| \left\|
u^{-1}\right\| =c$. It may be assumed that $ux_{\alpha }\in Y$ is
represented as a block 
\begin{equation*}
ux_{\alpha }=\sum\nolimits_{k=1}^{n(\alpha )}a_{k}^{\alpha }y_{\beta
_{k}\left( \alpha \right) }
\end{equation*}
for some sequence of rational scalars $\left( a_{k}^{\alpha }\right) _{k\leq
n(\alpha )}$ and finite $n(\alpha )$. Moreover, it may be assumed that this
blocks are not intersected, i.e. that any member of a given sequence $\left(
\beta _{k}\left( \alpha \right) \right) _{k\leq n(\alpha )}$ belongs only to
this block.

Since $\varkappa $ is uncountable, among such blocks there is an infinite
number of identical ones, that differs only in sequences $\left( \beta
_{k}\left( \alpha \right) \right) _{k\leq n(\alpha )}$. Let $A\subset
\varkappa $ be a such that all elements $\{ux_{\alpha }:\alpha \in A\}$ are
represented by those identical blocks: 
\begin{equation*}
ux_{\alpha }=\sum\nolimits_{k=1}^{n}a_{k}y_{\beta _{k}\left( \alpha \right) }%
\text{ \ for }\alpha \in A.
\end{equation*}

Let $\left( b_{i}\right) $ be a sequence of scalars; $a=\underset{1\leq
k\leq n}{\max }\left( \left| a_{k}\right| \right) $.

Then for any finite subset $A^{\prime }\subset A$, $A^{\prime }=\{\alpha
_{i}\}_{i=1}^{m}$, because of unconditionality of sequences $\left(
x_{\alpha }\right) $ and $\left( y_{\alpha }\right) $, 
\begin{eqnarray*}
\left\| \sum\nolimits_{i=1}^{m}b_{i}x_{\alpha _{i}}\right\| &\geq &c\left\|
\sum\nolimits_{i=1}^{m}b_{i}\left( \sum\nolimits_{k=1}^{n}a_{k}y_{\beta
_{k}\left( i\right) }\right) \right\| \\
&\geq &c\left\| \sum\nolimits_{i=1}^{m}b_{i}nay_{\beta _{1}\left( i\right)
}\right\| \geq cna\left\| \sum\nolimits_{i=1}^{m}b_{i}y_{\beta _{1}\left(
i\right) }\right\| .
\end{eqnarray*}

Analogously, the converse inequality that proves the theorem may be obtained.
\end{proof}

\section{Further constructions of spaces with the Tsirelson properties}

It was noted that superstable classes are only a constituent part of all
classes of finite equivalence. There exist classes of finite equivalence, in
which every space fails to have the Tsirelson property, and which, at the
same time, are not superstable.

\begin{example}
Let $2<p<q<\infty $. Consider a sequence of Banach spaces 
\begin{equation*}
Y_{1}=l_{2};\text{ }Y_{2}=l_{p}\left( Y_{1}\right) ;\text{ }%
Y_{3}=l_{2}\left( Y_{2}\right) ;\text{ }...;\text{ }Y_{2n}=l_{p}\left(
Y_{2n-1}\right) ;\text{ }Y_{2n+1}=l_{2}\left( Y_{2n}\right) ;\text{ }...
\end{equation*}
Similarly to the theorem 3 it may be shown that for a limit space $Y_{\infty
}=\underset{\rightarrow }{\lim }Y_{n}$ its\ index of divisibility $%
Index(Y_{\infty })=\{2,p\}$. According to the theorem 2 there exists a space 
$Z\sim _{f}Y_{\infty }$ that has the Tsirelson property. Hence, $Y_{\infty }$
is not isomorphic to a stable space.

Consider a space $X_{pq}=Y_{\infty }\oplus _{q}l_{q}$. Certainly, this space
is not superstable. However, every space $Z$ from the class $\left(
X_{pq}\right) ^{f}$ contains a subspace, isomorphic to $l_{r}$. Indeed, if $%
Z\in \left( X_{pq}\right) ^{f}$ then $Z$ is of kind $Z=Z_{1}\oplus _{q}Z_{2}$%
, where $Z_{1}\in \left( Y_{\infty }\right) ^{f}$; $Z_{2}\in \left(
l_{q}\right) ^{f}$ (since $S(Z)=\{2,p,q\}$ and $q\notin S(Y_{\infty })$.
Obviously, $Z_{2}$ contains a subspace, isomorphic to $l_{q}$.
\end{example}

\begin{theorem}
Let $X$ be a Banach space. If the corresponding class $X^{F}$ of crudely
finite equivalence is not crudely superstable then there exists a Banach
space $Y$, which is finite representable in $X$ and has the Tsirelson
property.
\end{theorem}

\begin{proof}
According to the proof of theorem 2, it is enough to check a such space $%
Z<_{f}X$ that the corresponding class $Z^{f}$ contains two spaces: $Z_{1}$
and $Z_{2}$ with $T\left( Z_{1}\right) \cap T\left( Z_{2}\right)
=\varnothing $. Recall that $T(Z)$ denotes a set of all $p\in \left[
1,\infty \right] $ such that $Z$ contains a subspace isomorphic to $l_{p}$.

Since $X^{F}$ is not crudely superstable, there exists $\left( W,\left(
w_{n}\right) \right) \in IS(X)$ such that $\frak{w}=\left( w_{n}\right) $ is
a spreading sequence, non-equivalent to any symmetric one.

Consider a tower $\left\lfloor \frak{w}\right\rfloor $ and all separable
spaces from it. By the results of previous section, there are an uncountable
set of such spaces that are pairwice non-isomorphic.

Let $\{\frak{w}_{\sigma }:\sigma \in \Sigma _{0}\}$ be their numeration,
where $\Sigma _{0}\subset \Sigma $ is a set of all rearrangements of $\omega 
$, that are not almost identical.

It will be assumed that $\frak{w}_{\sigma }=\{w_{\sigma n}:n<\omega \}$
forms a spreading basis of a space $W_{\sigma }$.

It will be shown that $\cap \{T\left( \frak{w}_{\sigma }\right) :\sigma \in
\Sigma _{0}\}=\varnothing $. Since every $W_{\sigma }<_{f}X$, this will
prove the theorem.

Suppose that this intersection is not empty, i.e. that there exists such $%
p\in \left[ 1,\infty \right] $ that $l_{p}$ is isomorphic to a subspace of
every $W_{\sigma }$ ($\sigma \in \Sigma _{0}$).

The natural basis $\left( e_{n}\right) $ of $l_{p}$ is reproducible in the
terminology of [11], i.e. if $l_{p}$ is isomorphic to a subspace of a Banach
space $Z$ with a basis $\left( z_{n}\right) $ then there exists a such
isomorphical embedding $u:l_{p}\rightarrow Z$ that 
\begin{equation*}
u\left( e_{k}\right) =\sum\nolimits_{i=m_{k}+1}^{m_{k+1}}a_{i}z_{i}
\end{equation*}
for a sequence of scalars $\left( a_{i}\right) $ and a sequence of naturals $%
m_{1}<...<m_{k}<...$.

Hence there exists a such isomorphic embedding $u_{\sigma }:l_{p}\rightarrow
W_{\sigma }$ that 
\begin{equation*}
u_{\sigma }\left( e_{k}\right) =\sum\nolimits_{i=m_{k}\left( \sigma \right)
+1}^{m_{k+1}\left( \sigma \right) }a_{i}^{\left( \sigma \right) }w_{\sigma
i}.
\end{equation*}

Because of $card\left( \Sigma _{0}\right) >\omega $, and $\{\frak{w}_{\sigma
}:\sigma \in \Sigma _{0}\}$ are spreading sequences, it may be chosen an
uncountable subset $\Sigma _{0}^{\left( 1\right) }\subset \Sigma _{0}$ such
that for all $\sigma \in \Sigma _{0}^{\left( 1\right) }$%
\begin{equation*}
u_{\sigma }\left( e_{1}\right) =\sum\nolimits_{i=1}^{m_{1}}a_{i}w_{\sigma i},
\end{equation*}
where $m_{1}$ and $\left( a_{i}\right) _{i=1}^{m_{1}}$ does not depend on $%
\sigma \in \Sigma _{0}^{\left( 1\right) }$.

Proceeding by induction, it may be chosen an uncountable subset $\Sigma
_{0}^{\left( 2\right) }\subset \Sigma _{0}^{\left( 1\right) }$ such that for
all $\sigma \in \Sigma _{0}^{\left( 2\right) }$ 
\begin{equation*}
u_{\sigma }\left( e_{k}\right)
=\sum\nolimits_{i=m_{1}+1}^{m_{2}}a_{i}w_{\sigma i},
\end{equation*}
etc. As a result it will be obtained a sequence $\Sigma _{0}^{\left(
1\right) }\supset ...\supset \Sigma _{0}^{\left( n\right) }\supset ...$ of
uncountable subsets of $\Sigma $, whose intersection is not empty (since $%
\omega $ is not confinal with $\omega _{1}$) such that for all $\sigma \in
\Sigma _{0}^{\infty }=\cap _{n}\Sigma _{0}^{\left( n\right) }$ 
\begin{equation*}
u_{\sigma }\left( e_{k}\right)
=\sum\nolimits_{i=m_{k}+1}^{m_{k+1}}a_{i}w_{\sigma i},
\end{equation*}
where $\left( n_{k}\right) _{k<\omega }$ and $\left( a_{k}\right) _{k<\omega
}$ does not depend on $\sigma \in \Sigma _{0}^{\infty }$.

It was assumed that all sequences $\{\{u_{\sigma }\left( e_{k}\right)
:k<\omega \}:\sigma \in \Sigma _{0}^{\infty }\}$ are equivalent to the
natural basis $\left( e_{n}\right) $ of $l_{p}$. Hence there exists a such
constant $c\in \left( 0,\infty \right) $ that for every finite sequence $%
\left( \xi _{j}\right) _{j=1}^{n}\ $of scalars and every $\sigma \in \Sigma
_{0}^{\infty }$ 
\begin{eqnarray*}
c^{-1}\left\| \sum\nolimits_{k=1}^{n}\xi _{k}\left(
\sum\nolimits_{i=m_{k}+1}^{m_{k+1}}a_{i}w_{\sigma i}\right) \right\| &\leq
&\left\| \sum\nolimits_{k=1}^{n}\xi _{k}\left(
\sum\nolimits_{i=m_{k}+1}^{m_{k+1}}a_{i}w_{i}\right) \right\| \\
&\leq &c\left\| \sum\nolimits_{k=1}^{n}\xi _{k}\left(
\sum\nolimits_{i=m_{k}+1}^{m_{k+1}}a_{i}w_{\sigma i}\right) \right\|
\end{eqnarray*}

Certainly, this implies that the sequence $\{w_{i}:i<\omega \}$ is
equivalent to a symmetric sequence.

Indeed, assume that it is not equivalent to any symmetric sequence. Then
there exists a double sequence $\{\{\xi _{k}^{\left( n\right) }:k\leq
m_{n}\}:n<\infty \}$ such that for all $n<\infty $ $\left\|
\sum\nolimits_{k=1}^{m_{n}}\xi _{k}^{\left( n\right) }w_{k}\right\| =1$ and $%
\underset{n\rightarrow \infty }{\lim }\left\| \sum\nolimits_{k=1}^{m_{n}}\xi
_{k}^{\left( n\right) }w_{\sigma k}\right\| =0$ for some rearrangement $%
\sigma $.

The expression $\left\| \sum\nolimits_{k=1}^{n}\xi _{k}\left(
\sum\nolimits_{i=m_{k}+1}^{m_{k+1}}a_{i}w_{i}\right) \right\| $ may be
presented in a form 
\begin{equation*}
\left\| \sum\nolimits_{k=1}^{n}\xi _{k}\left(
\sum\nolimits_{i=m_{k}+1}^{m_{k+1}}a_{i}w_{i}\right) \right\| =\left\|
\sum\nolimits_{k=1}^{n}\xi _{k}a_{m_{k}}w_{k}+z\right\| .
\end{equation*}
Similarly, 
\begin{equation*}
\left\| \sum\nolimits_{k=1}^{n}\xi _{k}\left(
\sum\nolimits_{i=m_{k}+1}^{m_{k+1}}a_{i}w_{\sigma i}\right) \right\|
=\left\| \sum\nolimits_{k=1}^{n}\xi _{k}a_{m_{k}}w_{\sigma k}+z^{\prime
}\right\| .
\end{equation*}

Let $\xi _{k}a_{m_{k}}=\xi _{k}^{\left( n\right) }$. This is possible since $%
\left( \xi _{k}\right) $ is an arbitrary finite sequence of scalars. Then 
\begin{eqnarray*}
c^{-1}\left\| \sum\nolimits_{k=1}^{n}\xi _{k}a_{m_{k}}w_{\sigma k}+z\right\|
&\leq &\left\| \sum\nolimits_{k=1}^{n}\xi _{k}a_{m_{k}}w_{k}+z\right\| \\
&\leq &c\left\| \sum\nolimits_{k=1}^{n}\xi _{k}a_{m_{k}}w_{\sigma
k}+z\right\| ,
\end{eqnarray*}
where $z$ in both cases may be assume to be the same element. So, 
\begin{equation*}
\frac{c^{-1}\left\| \sum\nolimits_{k=1}^{n}\xi _{k}^{\left( n\right)
}w_{\sigma k}+z\right\| }{\left\| \sum\nolimits_{k=1}^{m_{n}}\xi
_{k}^{\left( n\right) }w_{k}\right\| }\leq \frac{\left\|
\sum\nolimits_{k=1}^{n}\xi _{k}^{\left( n\right) }w_{k}+z\right\| }{\left\|
\sum\nolimits_{k=1}^{m_{n}}\xi _{k}^{\left( n\right) }w_{k}\right\| }
\end{equation*}
i.e. $c^{-1}\left\| e_{n}+z\right\| \leq \left\| e_{n}^{\prime }+z\right\| $%
, where $\left\| e_{n}^{\prime }\right\| =1$ and $\left\| e_{n}\right\| \ $%
infinitely increased. Certainly, this is impossible.
\end{proof}

\section{Ultracommutative Banach spaces}

Let $X$ be a Banach space; $D$, $E$ be ultrafilters. It will be said that
ultrapowers $\left( X\right) _{D}$ and $\left( X\right) _{E}$ are strongly
identical if there exists an isometry $i:\left( X\right) _{D}\rightarrow
\left( X\right) _{E}$ such that $i\circ d_{D}=d_{E}$, where $d_{D}$ and $%
d_{E}$ are canonical embeddings of $X$ into $\left( X\right) _{D}\ $and $%
\left( X\right) _{E}$ respectively. In a such case it will be written $%
\left\langle \left( X\right) _{D},d_{D}\right\rangle \equiv \left\langle
\left( X\right) _{E},d_{E}\right\rangle $.

\begin{definition}
A Banach space $X$ is said to be ultracommutative if for any pair of
ultrafilters $D$, $E$ over $\mathbb{N}$ ultrapowers $\left( X\right)
_{D\times E}$ and $\left( X\right) _{E\times D}$ are strongly identical.
\end{definition}

\begin{theorem}
Any ultracommutative Banach space $X$ is stable.
\end{theorem}

\begin{proof}
Let $\left( x_{n}\right) $ and $\left( y_{m}\right) $ be two sequences of
elements of $X$. The double sequence $z_{nm}=x_{n}+y_{m}$ generates a pair
of elements: $\left( z_{nm}\right) _{D\left( n\right) \times E\left(
m\right) }\in \left( X\right) _{D\times E}$ and $\left( z_{nm}\right)
_{E\left( m\right) \times D\left( n\right) }\in \left( X\right) _{E\times D}$%
. Since $X$ is ultracommutative, these elements are of equal norms. Indeed,
let $i:\left( X\right) _{D\times E}\rightarrow \left( X\right) _{E\times D}$
be an isometry. then, by the ultracommutativity, $iz_{nm}=z_{nm}$ and $%
i\left( z_{nm}\right) _{D\left( n\right) \times E\left( m\right) }=\left(
z_{nm}\right) _{E\left( m\right) \times D\left( n\right) }$. Certainly, this
implies that $\lim_{D\left( n\right) }\lim_{E\left( m\right) }\left\|
x_{n}+y_{m}\right\| =\lim_{E\left( m\right) }\lim_{D\left( n\right) }\left\|
x_{n}+y_{m}\right\| $.
\end{proof}

\begin{corollary}
In a general case iterated ultrapowers $\left( X\right) _{D\times E}\ $and $%
\left( X\right) _{E\times D}$ are not strongly identical.
\end{corollary}

The converse is not true.

\begin{theorem}
Let $X$ be an ultracommutative Banach space. Then each its spreading model
is isometric to a some $l_{p}$ ($1\leq p<\infty $).
\end{theorem}

\begin{proof}
As was noted in the proof of theorem 6 (cf. also [12]) every spreading model 
$sm(X,\left( x_{n}\right) ,D)$ of $X$ may be obtained by using iterated
ultrapowers; its natural basis is constructed by induction: 
\begin{equation*}
e_{1}=\left( x_{n}\right) _{D}\in \left( X\right) _{D};\text{ \ }%
e_{2}=\left( x_{n}\right) _{D\times D}\in \left( \left( X\right) _{D}\right)
_{D},\text{ \ etc.}
\end{equation*}

Assume that $X$ is ultracommutative. Then, by the preceding theorem it is
stable, in particular every its spreading model has a symmetric basis.

For $e_{1}$, $e_{2}$ as above and scalars $a$ and $b$ put: 
\begin{equation*}
f\left( a,b\right) =\left\| ae_{1}+be_{2}\right\| .
\end{equation*}

Certainly, the function $f$ is

homogeneous: $f\left( \lambda a,\lambda b\right) =\lambda f\left( a,b\right) 
$;

symmetric: $f\left( a,b\right) =f\left( b,a\right) $;

monotone: $f\left( a,b\right) \leq f\left( c,c\right) $ provided $a<c$ and $%
b<d$;

satisfies the norming condition $f(0,1)=1$ and, at least, from the
ultracommutativity follows: 
\begin{equation*}
f\left( a,f\left( b,c\right) \right) =f\left( f\left( a,b\right) ,c\right) .
\end{equation*}

Hence, according to the Kolmogoroff-Nagumo theorem (cf. [13] and [14]) or
from the same result, obtained independently by Bohnenblust [15], either
there exists such $p<\infty $ that $f\left( a,b\right) =\left(
a^{p}+b^{p}\right) ^{1/p}$ or $f\left( a,b\right) =\max \left( a,b\right) $.
Since $X$ is stable, the last case is impossible. So, $sm(X,\left(
x_{n}\right) ,D)=l_{p}$.
\end{proof}

\begin{corollary}
Spaces $L_{p}$ ($1\leq p\neq 2\leq \infty $) are not ultracommutative.
\end{corollary}

\begin{proof}
According to [16] for $1<p<2$ every $L_{p}$ contains a subspace $Y_{p}$ with
a symmetric basis that is complementably universal in the class of all
subspaces of $L_{p}$ with an unconditional basis. Certainly, $Y_{p}$ is not
isomorphic to any $l_{r}$. For $p>2$ the result follows by duality.

Indeed, for superreflexive Banach spaces $\left( X^{\ast }\right)
_{D}=\left( \left( X\right) _{D}\right) ^{\ast }$, that implies that $X$ and 
$X^{\ast }$ in a superreflexive case either both are non-ultracommutative,
or both enjoy this property.

To close the proof notice that $L_{1}$ for every $p\in \left( 1,2\right) $
contains a subspace, isometric to $L_{p}$.
\end{proof}

\begin{remark}
In a general case even if $X$ is ultracommutative, its ultrapower $\left(
X\right) _{D}$ (by a countably incomplete ultrafilter) does not have this
property: if its $l_{p}$-spectrum contains a point $p\neq 2$ then $\left(
X\right) _{D}$ contains a subspace, isometric to $L_{p}$. However, there are
ultracommutative spaces $X$ with ultracommutative ultrapowers that are not
isomorphic to $L_{2}$.
\end{remark}

\begin{example}
Consider a space $Z_{\frak{r}}=\left( \sum \oplus l_{p_{i}}^{\left(
n_{i}\right) }\right) _{2}$, where $p_{i}\rightarrow 2$; $n_{i}\rightarrow
\infty $; $\left| 2-p_{i}\right| \log n_{i}\rightarrow \infty $. Certainly, $%
Z_{\frak{r}}$ is stable; every space $X$, which is finitely equivalent to $%
Z_{\frak{r}}$ is of kind $X=Z_{\frak{r}}\oplus _{2}l_{2}\left( dim\left(
X\right) \right) $. In particular, all separable spaces from the class $%
\left( Z_{\frak{r}}\right) ^{f}$ are pairwice isometric. It is of interest
to point out that the space $Z_{\frak{r}}$ is stable (and also is
superstable) by the cardinal criterion (theorem 11); $card\left( \left( Z_{%
\frak{r}}\right) ^{f}\cap \mathcal{B}_{\varkappa }\right) =1<2^{\varkappa }$
for all cardinals $\varkappa $.

Notice, that there are continuum of paired sequences $\frak{r}%
=\{p_{i},n_{i}\}_{i<\infty }$ that generate nonisomorphic spaces of kind $Z_{%
\frak{r}}$.
\end{example}

\begin{example}
It is easy to check that any space of kind $\left( \sum \oplus A_{k}\right)
_{p}$, where $\left( A_{k}\right) $ is a sequence of finite dimensional
Banach spaces and $1\leq p<\infty $, is ultracommutative.\ 
\end{example}

\section{References}

\begin{enumerate}
\item  Tsirelson B.S. \textit{Not every Banach space contains }$l_{p}$%
\textit{\ or }$c_{0}$ (in Russian)\textit{,} Funct. Analysis and Appl. 
\textbf{8:2} (1974) 57-60

\item  Tokarev E.V. \textit{Structural theory of Banach spaces,} Ukrainian
Mathematical Congress-2001, International conference on functional analysis
(22-26August 2001 Kiev) Book of abstracts, Kiev (2001) 96-97

\item  Krivine J.-L., Maurey B. \textit{Espaces de Banach stables}, Israel
J. Math. \textbf{39} (1981) 273-295

\item  Raynaud Y. \textit{Espaces de Banach superstables}, C. R. Acad. Sci.
Paris, S\'{e}r. A. \textbf{292:14} (1981) 671-673

\item  Brunel A., Sucheston L. \textit{On }$\mathit{B}$\textit{-convex
Banach spaces}, Math. System Theory \textbf{7} (1973) 294-299

\item  Heinrich S. \textit{Ultraproducts in Banach space theory}, J. Reine
Angew. Math. \textbf{313} (1980) 72-104

\item  Pe\l czy\'{n}ski A. \textit{Projections in certain Banach spaces},
Studia Math. \textbf{19} (1960) 209 - 228

\item  Tokarev E.V. \textit{Notion of indiscernibless and on the }$\mathit{p}
$\textit{-Mazur property in Banach spaces} (transl. from Russian), Ukrainian
Mathematical Journal \textbf{37} (1985) 180-184

\item  Pe\l czy\'{n}ski A. \textit{Universal bases}, Studia Math. \textbf{%
32:3} (1969) 247 - 268

\item  Lindenstrauss J., Tzafriri L. \textit{Classical Banach spaces},
Lecture Notes in Math. \textbf{338 }(1973)

\item  Lindenstrauss J., Pe\l czy\'{n}ski A. \textit{Contribution to the
theory of the classical Banach spaces}, J. Funct. Anal. \textbf{8} (1971)
225 - 249

\item  Beauzamy B., Laprest\'{e} J.T. \textit{Mod\`{e}les \'{e}tal\'{e}s et
espaces de Banach}, Public. des D\'{e}part. Math\'{e}m. Univ. Claude
Bernard, Lyon-1, \textbf{4A} (1983) 1-199

\item  Kolmogoroff A. \textit{Sur notion de la moyene}, Rend. Acad. dei
Lincei \textbf{12:6} (1930) 388-391

\item  Nagumo M. \textit{Uber eine Klasse des Mittelwerte}, Japan J. Math. 
\textbf{7 }(1930)

\item  Bohnenblust H.F. \textit{An axiomatic characterization of }$L_{p}$%
\textit{\ - spaces}, Duke Math. J. \textbf{6} (1940) 627-640
\end{enumerate}

\end{document}